\documentclass [12pt]{article}

\usepackage {hyperref}
\newcommand{\tab}{\makebox[4em]{}}

\begin{document}
\begin{center}
\textbf{Infinite sequences in the framework of classical logic  }
\end{center}

\bigskip

\textbf{V.V. Ivanov}

N.N. Semenov Institute of Chemical Physics, Russian Academy of Sciences,
Kosygin str.4, Moscow, 119991, Russia (e-mail:
\href{mailto:vvi007@gmail.com}{vvi007@gmail.com}; vvivanov@chph.ras.ru)

\bigskip

\textbf{Abstract}

Infinite sequences are considered in the framework of classical logic from a
new point of view.\tab


\bigskip

\textbf{Introduction }

\bigskip

The problem of the existence of infinite sequences in the class of logical
objects is considered as solved positively and is not discussed in the
literature. Although the existence of infinite sequences is obviously
unprovable, non-existence of them is also obviously unprovable and this is
considered as a base for studying of properties of infinite sequences.
However more detailed analysis of the given problem in the framework of
classical logic shows that exactly this unprovability of both contradictory
statements is the final element of the proof that infinite sequences exist
only in the class of illogical objects.

\bigskip

\textbf{Initial states }

\bigskip

(1) Mathematical object (further in abbreviated form, object) is the subject
matter, which has provable properties (at least a part of properties is
provable). Objects are considered in the framework of classical logic.

\bigskip

(2) Class of objects is the object, which consists of other objects. All objects
of one class have some class of common properties. Properties are such
objects, which are indications of belonging (or not belonging) of other
objects to the considered classes. The logic is the method of consideration,
which allows dividing objects into classes. The proofs are objects, which
determine by logical method belonging (or not belonging) of the considered
objects to the considered classes. The object can only belong or not belong
to any class; the classical logic therefore is used.

\bigskip

(3) The belonging of object to class is also referred to as \textit{the existence} of object in the
given class. The class, in which the existence of object is considered, must
be defined in an explicit or implicit form. The existence of object is
ambiguous without this definition. According to the law of contradiction,
object can only exist or not exist in the same class. Aristotle's law of
contradiction [1] is strictly expressed in this formulation.

\bigskip

(4) Objects are divided into classes of \textit{logical} and \textit{illogical} objects. The object is referred to
as illogical, if it has \textit{unprovable} properties. For unprovable property of object, it
is unprovable that it exists in the class of properties of this object and
is unprovable that it does not exist in the same class (there is no
sufficient condition for considering that object has the given unprovable
property and no sufficient condition for considering that object does not
have this property). It is clear that the law of contradiction is
inapplicable to the pair \textit{the existence }-- \textit{non-existence} of unprovable property in class of properties
of illogical object. Indeed, one cannot consider that the law of
contradiction holds, since a sufficient condition for the choice of only the
existence or only non-existence is absent. One also cannot consider that the
law of contradiction is violated, since, because of unprovability of the
existence and non-existence, it is impossible to prove that one follows from
another.

\bigskip

(5) The belonging to the class of illogical objects is provable property of
object. If \textit{this belonging} is unprovable, the object has unprovable property, and this
proves its belonging to the class of illogical objects, hence \textit{this belonging} is provable.
However, logical inconsistency of provability in classical logic is
impossible and this proves the initial statement. Self-evident impossibility
of the choice of the existence or non-existence of any property of the
object in the class of its properties serves the proof in this case.

\bigskip

(6) Illogical objects cannot have other provable properties, besides belonging to the class of illogical objects. It is clear that properties, for which direct proof of their unprovability exists (under assumption that the rest of the properties are provable), are unprovable. The remaining properties are also unprovable. Indeed, even if there are provable properties of logical objects among the remaining properties, it does not follow from anything that they are also provable, when they belong to illogical objects. Thus, the remaining properties have unprovable property (belonging to the class of provable properties), consequently, they are unprovable. It is clear that logical analysis of illogical objects is not unreasonable only for proving that these objects belong to the class of illogical objects.

\bigskip

\textbf{Examination of infinite sequences }

\bigskip

\textbf{Theorem 1.} \textit{ Infinite periodic sequences belong to the class of illogical objects. }

\bigskip

\textbf{Proof.} Let us assume that infinite periodic sequence $B = p,\,p,\,...$ ($p$
is period of this sequence) is a logical object. Then properties of $B$ can
be analyzed logically. Let $p_1 $ represent the first period of $B$. Object
$p_1 $ is not property of $B$, but object, which exists in class of objects
of $B$. Let $E\left( {p_1 } \right)$ represent the property of $B$, which
consists in the existence of the first period $p_1 $ in class of periods of
$B$. Then $\neg E\left( {p_1 } \right)$ represent non-existence of the first
period $p_1 $ in class of periods of $B$. It is clear that the existence of
$E\left( {p_1 } \right)$ in class of periods of $B$ (the existence of the
existence of $p_1 $ in class of periods of $B)$ is nothing but $E\left( {p_1
} \right)$. It is also clear that non-existence of $E\left( {p_1 } \right)$
in class of periods of $B$ (non-existence of the existence of $p_1 $ in
class of periods of $B)$ is nothing but $\neg E\left( {p_1 } \right)$. The
unprovability of $E\left( {p_1 } \right)$ and unprovability of $\neg E\left(
{p_1 } \right)$ is clear, because although the first period is written in $B
= p,\,p,\,...$, nothing forbids to think that this period was excluded.
Since object $B$ has unprovable property $E\left( {p_1 } \right)$, initial
assumption is incorrect, consequently, $B$ belongs to the class of illogical
objects.

\bigskip

\textbf{Corollary. }If some infinite sequence exists in the class of logical objects,
then Cantor's diagonal method [2] is inapplicable in general case.

\bigskip

\textbf{Theorem 2.} \textit{ Infinite periodic sequences belong to the class of illogical objects. }

\bigskip

\textbf{Proof.} Let us consider infinite sequence $Q = a,b,c,...$, for which
algorithm $s$ is known, which determines any subsequent term, if the
previous one is known. Let us assume that $Q$ is logical object, then
properties of $Q$ can be analyzed logically. We denote $S$ class of all
terms, which can be obtained by algorithm $s$. We denote $E\left( S \right)$
the property of $Q$, which consists in the existence of $S$ in class of
objects of $Q$. Then $\neg E\left( S \right)$ represent non-existence of $S$
in class of objects of $Q$. It is clear that the existence of property
$E\left( S \right)$ in class of properties of $Q$ (the existence of the
existence of object $S$ in class of objects of $Q)$ is nothing but $E\left(
S \right)$. It is also clear that non-existence of property $E\left( {p_1 }
\right)$ in class of properties of $Q$ (non-existence of the existence of
object $S$ in class of objects of $Q)$ is nothing but $\neg E\left( S
\right)$. The unprovability of $E\left( S \right)$ and unprovability of
$\neg E\left( S \right)$ is clear, because it is possible considering that
sequence $Q$ contains all terms of $S$ (because any term of $S$ can be
transferred to $Q)$, however, it is also possible considering that sequence
$Q$ contains not all terms of $S$ (because transferring of terms of $S$ to
$Q$ cannot be completed). Since object $Q$ has unprovable property $E\left(
S \right)$, initial assumption is incorrect, consequently, $Q$ belongs to
the class of illogical objects.

\bigskip

\textbf{Conclusions}

\bigskip

The paradoxes concerned with objects, reveal unprovable properties of these
objects and prove that these objects belong to the class of illogical
objects. For example, should we consider as rational numbers all positive
numbers, which are expressed by vulgar fraction $\frac{m}{n}$? If no largest
$m = m_{\lim } $ exists in the class of natural numbers, then rational
numbers can become irrational. However if $m_{\lim } $ exists in the class
of natural numbers, then there exist fractions $\frac{m}{n}$, which are not
rational numbers. A sufficient condition for the choice of the existence or
non-existence of $m_{\lim } $ is absent. Thus, there is unprovable property
of positive rational numbers (in class of properties of these numbers),
which consists in the existence of $m_{\lim } $ (in the class of natural
numbers). Hence, positive rational numbers are the objects, which belong to
the class of illogical objects. Well-known paradoxes of Cantor, Russell and
other paradoxes of discrete infinity can be analyzed along the same scheme,
and each of them can be an example of proof that considered special case of
discrete infinity belongs to the class of illogical objects. It is obvious,
that the given scheme holds not only for paradoxes of discrete infinity. For
example, the statement "I am lying" has property 1 (in the class of Boolean
values 0 and 1), the existence of which (in class of properties of the
statement "I am lying") cannot be proved or refuted. The existence of
unprovable property means, that the considered statement is illogical
object.

\bigskip

The proofs of the present paper are not based on properties of numbers,
therefore incompleteness of axioms of arithmetic cannot be counterargument
for them. On the contrary, Godel's incompleteness theorem [3] demands the
existence of nondenumerability in the class of logical objects, while
nondenumerability is property of infinite sets, which belong to the class of
illogical objects.

\bigskip

The next author's paper will be devoted to the proof that in contrast to infinite sequences discrete infinity can exist in the form of infinite numerical matrices in the class of \textit{logical} objects.

\bigskip

\textbf{References\tab }

1. Aristotle, \textit{Metaphysics}, Book IV, part III.

2. Cantor, G.: \"{U}ber eine elementare Frage der Mannigfaltigkeitslehre.
Jahresbericht der Deutschen Math. Vereinigung \textbf{I }, 75-78
(1890-1891).

3. Goedel, K.: \"{U}ber
formal unentscheidbare S\"{a}tze der Principia Mathematica und verwandter
Systeme, I. Monatshefte f\"{u}r Mathematik und Physik \textbf{38}, 173-198
(1931).

\end{document}